\documentclass[11pt,leqno]{article}
\usepackage{amsmath}
\usepackage{amsthm}
\usepackage{amscd}
\usepackage{amssymb}
\usepackage[francais]{babel}
\usepackage[T1]{fontenc}
\usepackage{url}
\usepackage{enumerate}

\makeatletter
\renewcommand{\section}{%
 \@startsection{section}{1}{\z@}%
 {-3.5ex \@plus -1ex \@minus -.2ex}%
 {2.3ex \@plus.2ex}%
 {\centering\normalfont\Large\bfseries}}
\renewcommand{\@seccntformat}[1]{\csname the#1\endcsname.\hspace{.5em}}
\makeatother
\newtheoremstyle{standard}%  name
  {}%         Space above, empty = `usual value'
  {}%         Space below
  {\rmfamily\mdseries\itshape}% Body font
  {}%         Indent amount
  {\rmfamily\bfseries\upshape}% Thm head font
  {.\hspace{.5em}---}%        Punctuation after thm head
  {.6em}%     Space after thm head: " " = normal interword space;
  {\thmnumber{#2. }\thmname{#1}\thmnote{ \textmd{(#3)}}% Thm head spec
}
\newcounter{entry}[section]
\renewcommand{\theentry}{\thesection.\arabic{entry}}
\newcommand{\entry}[1][\hspace{-.5em}]{%
 \vspace{.5\baselineskip}\par%
 \refstepcounter{entry}%
 \noindent{\rmfamily\bfseries\upshape\theentry.\hspace{.5em}#1\hspace{.5em}---\hspace{.6em}}%
}
\theoremstyle{standard}
\newtheorem{theoreme}[entry]{Th\'eor\`eme}
\newtheorem{lemme}[entry]{Lemme}
\newtheorem{proposition}[entry]{Proposition}
\newtheorem{corollaire}[entry]{Corollaire}
\newtheorem{definition}[entry]{D\'efinition}
\newtheorem*{theoreme*}{Th\'eor\`eme}
\newtheorem*{prop-def*}{Proposition-d\'efinition}
\numberwithin{equation}{entry}

\title{M\'etriques de sous-quotient et th\'eor\`eme de Hilbert-Samuel arithm\'etique pour les faisceaux coh\'erents}
\author{Par \emph{Hugues Randriambololona}, \`a Paris.}
\date{}

\input{utile.def}

\begin{document}

\maketitle

\section*{Introduction}

Soient $K$ un corps de nombres, $\cO_K$ son anneau d'entiers, $\gX$ un
sch\'ema projectif sur $\spec\cO_K$
de fibre g\'en\'erique $\gX_K$ r\'eduite,
%et \'equidimensionnelle,
et $\overline{\cL}$ 
un $\cO_{\gX}$-module inversible ample dont les fibres
sur l'espace analytique r\'eduit $\gX(\C)$ sont,
de fa\c{c}on compatible \`a
la conjugaison complexe,
munies d'une m\'etrique continue \`a courbure semi-positive.

Si $\cC$ est un $\cO_{\gX}$-module coh\'erent de support
de dimension absolue $d\leq\dim\gX$, on note $[\cC]\in Z_d(\gX)$ le cycle
associ\'e \`a $\cC$, d\'efini comme suit~:
\begin{equation*}
[\cC]=\sum_{x\in\gX_{(d)}}(\lg_{\cO_{\gX,x}}\cC_x)[x].
\end{equation*}
La hauteur de ce cycle relativement \`a $\overline{\cL}$
peut \^etre d\'efinie comme le nombre
d'intersection arithm\'etique
$(\wc1(\overline{\cL})^d.[\cC])=\sum_{x\in\gX_{(d)}}(\lg_{\cO_{\gX,x}}\cC_x)(\wc1(\overline{\cL}_{|x})^d.[x])$,
suivant la
th\'eorie de l'intersection arithm\'etique de Gillet et Soul\'e
(cf.~\cite{GSAIT}, et \cite{Maillot} pour la
g\'en\'eralisation aux m\'etriques
non $C^\infty$).

Supposons maintenant donn\'es un $\cO_{\gX}$-module localement
libre hermitien $\overline{\cE}$, un sous-$\cO_{\gX}$-module
coh\'erent $\cF$ de $\cE$, et un morphisme
de $\cO_{\gX}$-modules $p:\cF\surj\cC$
qui fait de $\cC$
un quotient de $\cF$, ou encore, un <<sous-quotient>> de $\cE$.
Pour tout entier $n$,
le $\cO_K$-module $\Gamma(\gX,\cE\tens\cL^{\tens n})$
peut \^etre muni, au choix, des normes uniformes, ou bien,
si $\gX(\C)$ est muni d'une forme volume positive continue
compatible \`a la conjugaison complexe, des normes d'int\'egration $L^2$.
Par consid\'eration des normes restreintes, on obtient aussi
une structure de $\cO_K$-module norm\'e
sur $\Gamma(\gX,\cF\tens\cL^{\tens n})$, puis, par passage au quotient,
sur $\Gamma(\gX,\cC\tens\cL^{\tens n})$ si $n$ est assez grand.

Le r\'esultat principal de ce texte est la g\'en\'eralisation du
<<th\'eor\`eme de Hilbert-Samuel arithm\'etique>> de \cite{GSRR}
dans ce cadre des
faisceaux coh\'erents munis de <<normes de sous-quotient>>.
On montre en effet~:
\begin{theoreme*}
Sous les hypoth\`eses pr\'ec\'edentes, le degr\'e
d'Arakelov du $\cO_K$-module
norm\'e $\overline{\Gamma(\gX,\cC\tens\cL^{\tens n})}$ admet le
d\'eveloppement asymptotique
\begin{equation*}
\wdeg\overline{\Gamma(\gX,\cC\tens\cL^{\tens n})}=\frac{n^d}{d!}(\wc1(\overline{\cL})^d.[\cC])+o(n^d)
\end{equation*}
quand $n$ tend vers l'infini.
\end{theoreme*}

Ce th\'eor\`eme apporte une r\'eponse \`a deux questions
d\'ej\`a apparues dans la litt\'erature~:
\begin{enumerate}[(i)]
\item Dans \cite{GSRR} paragraphe 5.3.2, il est demand\'e si le
th\'eor\`eme de Hilbert-Samuel arithm\'etique prouv\'e pour un
faisceau localement libre pouvait s'\'etendre \`a n'importe quel
faisceau coh\'erent, et plus particuli\`erement \`a un sous-faisceau
coh\'erent d'un faisceau localement libre hermitien, muni des normes
restreintes~; ceci est bien un cas particulier du probl\`eme r\'esolu ici,
du moins en ce qui concerne le terme de degr\'e maximal de cette formule
de Hilbert-Samuel arithm\'etique.
\item Le faisceau structural $\cO_\Sigma$ d'un sous-sch\'ema ferm\'e $\Sigma$
de $\gX$ est naturellement muni d'une structure de quotient du faisceau
structural $\cO_\gX$ de $\gX$, ce dernier disposant d'une structure
m\'etrique \'evidente~; en appliquant le th\'eor\`eme pr\'ec\'edent
avec $\cC=\cO_\Sigma$ on obtient ainsi un
th\'eor\`eme de Hilbert-Samuel arithm\'etique pour les hauteurs de
sous-sch\'emas,
comme conjectur\'e dans \cite{these} A.2.1.
\end{enumerate}

Il serait int\'eressant d'examiner si le r\'esultat obtenu ici ne dispose
pas d'applications en th\'eorie de l'approximation diophantienne,
en autorisant par exemple l'utilisation de <<fonctions auxiliaires>>
tordues par des sections d'un faisceau coh\'erent
(non n\'ecessairement localement libre).

Au moyen d'un d\'evissage, on montre que pour prouver le th\'eor\`eme
en toute g\'en\'eralit\'e, il suffit de le faire dans le cas particulier
o\`u $\cC$ est de la forme $i_*\cM$, o\`u $i:Z\inj\gX$ est le
morphisme d'immersion d'un sous-sch\'ema ferm\'e int\`egre $Z$
de $\gX$ et $\cM$ un $\cO_Z$-module inversible. On se ram\`ene
alors \`a la version classique du th\'eor\`eme de Hilbert-Samuel
arithm\'etique, ou plut\^ot \`a sa g\'en\'eralisation au cas
singulier donn\'ee dans \cite{Zhang}, au moyen d'un r\'esultat de
comparaison de normes dont l'essentiel r\'eside en l'\'enonc\'e
de prolongement de sections holomorphes
d'un fibr\'e vectoriel hermitien sous-quotient avec
contr\^ole des normes qui constitue le point (ii) du th\'eor\`eme
suivant~:

\begin{theoreme*}
Soient $X$ un espace analytique $1$-convexe\footnote{par exemple,
compact, ou de Stein (cf. \S\S\ref{def_1-convexe}--\ref{ex_1-convexe})}
r\'eduit, $\overline{L}$ un $\cO_X$-module inversible hermitien
\`a courbure strictement positive, $\overline{E}$
un $\cO_X$-module localement libre hermitien de type fini, $i:Y\inj X$
un sous-espace
analytique ferm\'e r\'eduit de $X$, et $\overline{V}$ un $\cO_Y$-module
localement libre hermitien de type fini.
On suppose donn\'es un
sous-$\cO_X$-module
coh\'erent $F$ de $E$ et un morphisme surjectif
\begin{equation*}
p:F\surj i_*V
\end{equation*}
de $\cO_X$-modules coh\'erents.
Alors~:
\begin{enumerate}
\item
Pour
tout r\'eel $\epsilon>0$ et pour tout compact non
vide $B$ de $Y$, il existe un r\'eel $C>0$ et un compact
non vide $A$ de $X$ tels que
pour tout $n\geq 0$, pour tout $s\in\Gamma(Y,V\tens L^{\tens n})$
et pour tout $\widetilde{s}\in\Gamma(X,F\tens L^{\tens n})$
v\'erifiant $p(\widetilde{s})=i_*s$,
on ait
\begin{equation*}
\|\widetilde{s}\|_{L^\infty(A)}\geq C e^{-n\epsilon}\|s\|_{L^\infty(B)}.
\end{equation*}
\item
Il existe un entier $n_0$ et, pour
tout r\'eel $\epsilon>0$ et pour tout compact non
vide $A$ de $X$, un r\'eel $C'>0$ et un compact
non vide $B$ de $Y$, tels que
pour tout $n\geq n_0$ et pour tout $s\in\Gamma(Y,V\tens L^{\tens n})$
il existe $\widetilde{s}\in\Gamma(X,F\tens L^{\tens n})$
v\'erifiant $p(\widetilde{s})=i_*s$
et
\begin{equation*}
\|\widetilde{s}\|_{L^\infty(A)}\leq C'e^{n\epsilon}\|s\|_{L^\infty(B)}.
\end{equation*}
\end{enumerate}
\end{theoreme*}

On trouvera dans la litt\'erature des \'enonc\'es de prolongement
analogues \`a celui-ci (souvent dans le cas particulier $F=E$ et $V=E|_Y$)
\'etablis pour la plupart au moyen de la technique
des estim\'ees $L^2$ de H\"ormander, et pouvant \'eventuellement
donner un contr\^ole plus fort sur
les normes, mais au prix de certaines hypoth\`eses de
lissit\'e
(voir par exemple \cite{DemaillyOTM}, \cite{Manivel}, ou
encore \cite{these}, chap.~3) ou
de la connaissance \emph{a priori} de l'existence d'un
prolongement convenable de $s$ sur un
voisinage donn\'e de $Y$ (cf.~\cite{Zhang}, sect.~2).
L'originalit\'e du r\'esultat pr\'esent\'e ici est donc de
s'affranchir totalement de ces hypoth\`eses annexes.
Pour y parvenir on utilise la th\'eorie des espaces $1$-convexes
(apr\`es s'\^etre plac\'e sur le fibr\'e en disques du dual $L^\vee$)
et les techniques d'espaces de Fr\'echet
de la th\'eorie topologique des faisceaux coh\'erents sur
un espace analytique complexe, suivant en cela la m\'ethode introduite
par Bost dans \cite{Bost}, app.~A.

\vspace{.5\baselineskip}
\noindent\emph{Conventions}.\hspace{.5em}---
On utilisera ici la terminologie usuelle des espaces analytiques complexes
(non n\'ecessairement r\'eduits) et des fonctions continues, $C^\infty$,
holomorphes, plu\-ri\-sous\-har\-mo\-niques, etc. sur iceux, telle qu'elle est
rappel\'ee par exemple au d\'ebut de \cite{GrauertModif}
et de \cite{Nar} (on gardera notamment \`a l'esprit
le th\'eor\`eme 5.3.1 de \cite{FN}, selon lequel les
deux d\'efinitions raisonnables des fonctions
plu\-ri\-sous\-har\-mo\-niques co\"{\i}ncident).

Les m\'etriques hermitiennes sur les faisceaux analytiques
localement libres seront toujours suppos\'ees continues.
Un faisceau inversible hermitien $\overline{L}$
sur un espace analytique complexe
sera dit \`a courbure semi-positive (resp. strictement positive)
si, pour toute section $l$ de $L$ ne s'annulant pas sur un
ouvert $U$, la fonction $-\log\|l\|$ est \psh
(resp. strictement plu\-ri\-sous\-har\-mo\-nique) sur $U$.

\vspace{.5\baselineskip}
\noindent\emph{Remerciements}.\hspace{.5em}---
L'auteur remercie J.-B.~Bost pour l'int\'er\^et qu'il a
manifest\'e \`a l'\'egard de ce travail et pour la version
pr\'eliminaire de \cite{Bost} qu'il a bien voulu lui
communiquer.

\section{Sous-quotients}

\entry
\label{intro_sous-quotients}
Soit $\cA$ une cat\'egorie ab\'elienne. On rappelle que si $M$
est un objet de $\cA$, les sous-objets des quotients de $M$
s'identifient naturellement aux quotients des sous-objets de $M$,
ou encore aux gradu\'es $M_1/M_2$ des filtrations \`a deux termes
\begin{equation}
\label{filtration_2_termes}
M\supset M_1\supset M_2.
\end{equation}
De fa\c{c}on plus pr\'ecise, on appelle \emph{sous-quotient} de $M$ 
la donn\'ee d'une telle filtration \eqref{filtration_2_termes}.
On remarquera que les sous-quotients des objets de $\cA$ forment
une cat\'egorie additive (en g\'en\'erale non ab\'elienne), un
morphisme de $M\supset M_1\supset M_2$ dans $N\supset N_1\supset N_2$
\'etant un morphisme de $M$ dans $N$ qui envoie $M_1$ dans $N_1$
et $M_2$ dans $N_2$~; un tel morphisme induit alors naturellement
un morphisme entre les gradu\'es $M_1/M_2$ et $N_1/N_2$.

Par restriction, les sous-quotients d'un objet $M$ 
fix\'e forment aussi une cat\'egorie, l'ensemble des morphismes
de $M\supset M_1\supset M_2$ dans $M\supset M'_1\supset M'_2$ \'etant
non vide si et seulement si $M_1$ est inclus dans $M'_1$
et $M_2$ dans $M'_2$, et alors par
d\'efinition cet ensemble se r\'eduit au seul
morphisme d\'eduit de l'identit\'e de $M$. On prendra garde
toutefois qu'il est alors possible que le morphisme induit de $M_1/M_2$
dans $M'_1/M'_2$ soit un isomorphisme sans pour autant que le
morphisme de sous-quotients originel ne soit inversible~; une telle
situation se produit par exemple lorsque $M$
est un produit direct $A\times B$ et que l'on consid\`ere le
morphisme naturel entre les
sous-quotients $A\times B\supset A\times\{0\}\supset\{0\}$
et $A\times B\supset A\times B\supset \{0\}\times B$, les deux
gradu\'es associ\'es s'identifiant naturellement \`a $A$.

Il sera commode pour la suite d'utiliser la notion de sous-quotient
sous la forme suivante~:

\entry[D\'efinition.]
Soient $C$ et $E$ deux objets d'une cat\'egorie ab\'elienne.
On appelle \emph{structure de sous-quotient de $E$} sur $C$ la donn\'ee
d'un sous-quotient
\begin{equation}
E\supset E_1\supset E_2
\end{equation}
de $E$ et d'un isomorphisme
\begin{equation}
\phi:C\overset{\sim}{\longto}E_1/E_2
\end{equation}
de $C$ sur le gradu\'e associ\'e.

\entry
\label{caracterisations_sous-quotients}
Ainsi, munir $C$ d'une structure de sous-quotient de $E$ \'equivaut
encore \`a se donner, au choix~:
\begin{enumerate}[(i)]
\item un sous-objet $F$ de $E$ et un morphisme surjectif $p:F\surj C$  
(prendre $F=E_1$ et $p$ la surjection de noyau $E_2$
d\'eduite de $\phi^{-1}$)~;
\item un quotient $Q$ de $E$ et une injection $i:C\inj Q$
(prendre $Q=E/E_2$ et $i$ l'injection d\'eduite de $\phi$)~;
\item un complexe court
\begin{equation}
\cE:E'\overset{i}{\inj}E\overset{p}{\surj}E''
\end{equation}
avec $i$ injective, $p$ surjective, et $p\circ i=0$, et un isomorphisme
\begin{equation}
C\overset{\sim}{\longto}H(\cE)=\ker p/\im i
\end{equation}
de $C$ sur la cohomologie de $\cE$ (prendre $E'=E_2$ et $E''=E/E_1$).
\end{enumerate}
Compte tenu de cette derni\`ere caract\'erisation,
on dira aussi parfois que les sous-quotients
de $E$ sont les objets de cohomologie de $E$.

Enfin on prendra garde, comme cela a d\'ej\`a \'et\'e signal\'e \`a
la fin de \ref{intro_sous-quotients}, qu'il peut arriver
qu'un m\^eme objet $C$ admette
plusieurs structures diff\'erentes de sous-quotient
d'un m\^eme objet $E$.

\entry
\label{sous-quotients_emboites}
Soient $E$, $F$ et $G$ trois objets d'une cat\'egorie ab\'elienne.
On va montrer que la donn\'ee sur $F$ d'une structure de sous-quotient
de $E$ et sur $G$ d'une structure de sous-quotient de $F$ d\'etermine de
fa\c{c}on naturelle sur $G$ une structure de sous-quotient de $E$.

En effet, on dispose par hypoth\`ese d'un isomorphisme $F\simeq E_1/E_2$
o\`u $E_1\supset E_2$ sont deux sous-objets de $E$ et
d'un isomorphisme $G\simeq F_1/F_2$
o\`u $F_1\supset F_2$ sont deux sous-objets de $F$.
Par l'identification naturelle des sous-objets de $E_1/E_2$ aux
sous-objets de $E_1$ contenant $E_2$, $F_1$ et $F_2$ se rel\`event
canoniquement en deux sous-objets $\widehat{F_1}$ et $\widehat{F_2}$
de $E$ v\'erifiant les inclusions
\begin{equation}
\label{filtration_4_termes}
E_1\supset\widehat{F_1}\supset\widehat{F_2}\supset E_2,
\end{equation}
de sorte que $G$ peut bien \^etre muni de la structure de sous-quotient
de $E$ d\'efinie par la filtration \`a deux termes
$E\supset\widehat{F_1}\supset\widehat{F_2}$ et par
l'isomorphisme compos\'e $G\simeq F_1/F_2\simeq\widehat{F_1}/\widehat{F_2}$.

\entry[D\'efinition.]
La structure de sous-quotient de $E$ d\'efinie sur $G$ au paragraphe
pr\'ec\'edent sera appel\'ee \emph{structure de sous-quotient compos\'ee}
des structures de sous-quotient de $E$ sur $F$ et de $F$ sur $G$.

Dans cette situation, on dira aussi parfois que $G$ et $F$ sont
deux \emph{sous-quotients embo{\^\i}t\'es} de $E$.

\entry
Soient $(K,|.|)$ un corps valu\'e et $E$ un $K$-espace vectoriel
muni d'une norme $\|.\|$ compatible \`a $|.|$. Puisqu'un $K$-espace
vectoriel $V$ muni d'une structure de sous-quotient de $E$ peut \^etre
vu, au choix, comme un sous-espace d'un quotient ou comme un quotient
d'un sous-espace de $E$, il h\'erite naturellement d'une norme
obtenue \`a partir de la norme de $E$ en consid\'erant respectivement
la norme restreinte de la norme quotient ou la norme quotient de
la norme restreinte. On se convaincra facilement que ces deux
constructions donnent la m\^eme norme sur $V$~; de fa\c{c}on plus
pr\'ecise~:

\entry[D\'efinition.]
Soient $(E,\|.\|)$ un $(K,|.|)$-espace vectoriel norm\'e
et $V$ un $K$-espace vectoriel muni d'une structure de sous-quotient
de $E$ d\'efinie par une filtration $E\supset E_1\supset E_2$
et un isomorphisme $\phi:V\overset{\sim}{\longto}E_1/E_2$.
La norme de sous-quotient $\|.\|_{\mathrm{sq}}$ sur $V$ est
d\'efinie, pour $v\in V$, par la formule
\begin{equation}
\label{formule_norme}
\|v\|_{\mathrm{sq}}=\inf_{e\in E_1,\;\phi(v)=\overline{e}}\|e\|
\end{equation}
o\`u l'on a not\'e $\overline{e}$ la classe dans $E_1/E_2$
de l'\'el\'ement $e$ de $E_1$.

\'Etudions comment cette construction se comporte relativement
\`a la composition des structures de sous-quotient~:

\begin{proposition}
\label{transitivite_norme}
Soient $(E,\|.\|_E)$ un $(K,|.|)$-espace vectoriel
norm\'e, $F$ un $K$-espace vectoriel muni d'une structure de
sous-quotient de $E$, et $G$ un $K$-espace vectoriel muni d'une
structure de sous-quotient de $F$. Alors la norme sur $G$ sous-quotient
de la norme sur $F$ sous-quotient de $\|.\|_E$ co{\"\i}ncide
avec la norme sur $G$ sous-quotient de $\|.\|_E$ relativement
\`a la structure de sous-quotient compos\'ee.
\end{proposition}
On dira aussi que la formation des normes sous-quotient est
transitive dans les sous-quotients embo{\^\i}t\'es.
\begin{proof}
Notons $E\supset E_1\supset E_2$ et $F\supset F_1\supset F_2$
les filtrations et $p:E_1\surj F$ et $q:F_1\surj G$ les projections
de noyaux $E_2$ et $F_2$
d\'efinissant les structures de sous-quotient
consid\'er\'ees.
Notons aussi $\widehat{F_1}$ et $\widehat{F_2}$ les
sous-espaces de $E$ relevant $F_1$ et $F_2$ comme en
\eqref{filtration_4_termes}, de sorte que $\widehat{F_1}$
(resp.~$\widehat{F_2}$) est l'ensemble des \'el\'ements de $E$
dont l'image par $p$ appartient \`a $F_1$ (resp.~$F_2$)
et que l'application
compos\'ee $q\circ (p|_{\widehat{F_1}})$ de $\widehat{F_1}$
sur $G$ est bien une projection de noyau $\widehat{F_2}$
et d\'efinit la structure de sous-quotient compos\'ee.
Alors, par construction, 
notant $\|.\|_F$ la norme sur $F$ sous-quotient de $\|.\|_E$,
et $\|.\|_G$ la norme sur $G$ sous-quotient de $\|.\|_F$,
on a pour $g\in G$~:
\begin{equation}
\begin{split}
\|g\|_G&=\inf_{f\in F_1,\;q(f)=g}\|f\|_F\\
&=\inf_{f\in F_1,\;e\in E_1,\;q(f)=g,\;p(e)=f}\|e\|_E\\
&=\inf_{e\in \widehat{F_1},\;q(p(e))=g}\|e\|_E,
\end{split}
\end{equation}
ce qui correspond bien \`a la norme de $g$ pour la structure de
sous-quotient compos\'ee.
\end{proof}

Il sera aussi utile de pouvoir comparer les deux normes sous-quotient
de deux normes sur un m\^eme espace~:
\begin{proposition}
\label{compare_normes_sq}
Soient $C$ un r\'eel, $E$ un $K$-espace vectoriel
muni de deux normes $\|.\|$ et $\|.\|'$ v\'erifiant l'in\'egalit\'e
\begin{equation}
\label{tmp_compare_normes_sq}
\|.\|\leq C\|.\|'
\end{equation}
et $V$ un $K$-espace vectoriel
muni d'une structure de sous-quotient de $E$. Alors les normes
sous-quotient $\|.\|_{\mathrm{sq}}$ de $\|.\|$
et $\|.\|'_{\mathrm{sq}}$ de $\|.\|'$ sur $V$ v\'erifient
\begin{equation}
\|.\|_{\mathrm{sq}}\leq C\|.\|'_{\mathrm{sq}}.
\end{equation}
\end{proposition}
\begin{proof}
L'in\'egalit\'e \eqref{tmp_compare_normes_sq} passe \`a la borne
inf\'erieure dans \eqref{formule_norme}.
\end{proof}

\section{Prolongement de sections holomorphes avec contr\^ole des normes}

\entry\label{rappel_topologie_canonique}
On rappelle (cf.~\cite{GRStein}, chap.~V~\S6) que si $X$ est
un espace analytique complexe et $\cC$ un $\cO_X$-module coh\'erent,
il existe une unique topologie de Fr\'echet sur l'espace des sections
globales $\Gamma(X,\cC)$, appel\'ee topologie canonique,
qui rende continues toutes les applications
de localisation
\begin{equation}
\Gamma(X,\cC)\longto\cC_x
\end{equation}
en les fibres $x\in X$, lorsque le module $\cC_x$ sur l'alg\`ebre
analytique locale $\cO_{X,x}$ est muni de la topologie canonique
au sens de \cite{Jurchescu}
(\ie de la <<Folgentopologie>> de \cite{GRAS}).
On dispose en outre des propri\'et\'es suivantes~:
\begin{enumerate}[(i)]
\item Si $X$ est r\'eduit, et si $\overline{\cE}$ est un $\cO_X$-module
localement libre hermitien de type fini, la topologie canonique
sur $\Gamma(X,\cE)$ co{\"\i}ncide avec la topologie de la convergence
uniforme sur les compacts (cela r\'esulte par exemple de
\cite{GRStein}~V~\S6 th.~8, du caract\`ere local de la
condition d\'efinissant la topologie canonique, et des
compatibilit\'es \'evidentes aux sommes directes).
\item Si $\cC'\subset\cC$ est un sous-$\cO_X$-module coh\'erent, la
topologie canonique sur $\Gamma(X,\cC')$ est la topologie induite
par la topologie canonique sur $\Gamma(X,\cC)$ pour l'inclusion
naturelle (cela r\'esulte de l'assertion analogue pour la topologie
canonique sur les fibres~: \cite{GRAS}~II~\S2.7 Satz~9, pp.~97--99).
\item Si $i:Z\inj X$ est un sous-espace analytique ferm\'e
et si $\cD$ est un $\cO_Z$-module coh\'erent, l'identification
naturelle de $\Gamma(Z,\cD)$ et de $\Gamma(X,i_*\cD)$ est un
hom\'eomorphisme pour les topologies canoniques sur ces espaces
(c'est une cons\'equence directe des d\'efinitions).
\end{enumerate}

\entry
\label{def_1-convexe}
On utilisera aussi dans cette partie la th\'eorie des
espaces $1$-convexes \'elabor\'ee dans \cite{AG} et \cite{Nar}.
On suivra ici la pr\'esentation de \cite{VVT}~:
\begin{prop-def*}
Soit $X$ un espace analytique complexe. Les assertions
suivantes sont \'equivalentes~:
\begin{enumerate}[(i)]
\item Il existe sur $X$ une fonction d'exhaustion continue qui est
strictement \psh hors d'un certain compact. 
\item L'espace $X$ est holomorphiquement convexe et admet
un sous-ensemble analytique compact sans points isol\'es maximal.
\item L'espace $X$ est une modification propre
(au sens de \cite{GrauertModif}) d'un espace de Stein
en un nombre fini de points.
\item Pour tout faisceau analytique coh\'erent $F$ sur $X$
et tout entier $q\geq 1$,
le groupe de cohomologie $H^q(X,F)$ est de dimension finie.  
\end{enumerate}
Si l'une de ces conditions est v\'erifi\'ee, on dit
que $X$ est $1$-convexe. 
\end{prop-def*}
Si $X$ est $1$-convexe, l'espace de Stein $Y$ introduit en (iii)
est unique (\`a isomorphisme pr\`es)~; c'est le r\'eduit de
Cartan-Remmert de $X$.
Le sous-ensemble analytique compact maximal de (ii)
est alors le support du lieu
exceptionnel $S$ de la modification $X\longto Y$, et pour tout faisceau
analytique coh\'erent $F$ sur $X$ et tout entier $q\geq 1$ on
dispose d'un isomorphisme naturel $H^q(X,F)\simeq H^q(S,F|_S)$.
En outre, parmi les fonctions d'exhaustion continues dont l'existence
est assur\'ee par (i), il en est une qui est strictement \psh
en tout point du compl\'ementaire de $S$.

L'implication (i)$\Longrightarrow$(iv), qui est la seule dont nous aurons
r\'eellement besoin ici, est prouv\'ee dans \cite{AG}. Les autres r\'esultats
d\'ecoulent de \cite{Nar}, th.~V, et des %propri\'et\'es
caract\'erisations classiques des espaces de Stein.
%et de la th\'eorie des modifications.
On pourra consulter \cite{VVT} pour un expos\'e synth\'etique de tout ceci.

\entry\label{ex_1-convexe}
Parmi les exemples \'el\'ementaires d'espaces $1$-convexes on trouve
les espaces de Stein d'une part, et les espaces compacts d'autre part.
Tout sous-espace ferm\'e d'un espace $1$-convexe est $1$-convexe.

Un autre proc\'ed\'e de
construction d'espaces $1$-convexes dont nous aurons besoin est
celui qui suit.

\entry\label{notations_disque_1-convexe}
Pour tout  espace analytique complexe r\'eduit $X$ et
pour tout $\cO_X$-module inversible $L$
muni d'une m\'etrique continue $\|.\|$,
on note $V(X,L)$ l'espace total du fibr\'e en droites dual $L^\vee$
et, pour tout r\'eel $r>0$, $D_r(X,L)$ son fibr\'e en disques ouverts
de rayon $r$, relativement \`a la norme duale $\|.\|^\vee$.
Notons aussi
\begin{equation}
\pi:D_r(X,L)\surj X
\end{equation}
la projection naturelle
et 
\begin{equation}
\iota:X\inj D_r(X,L)
\end{equation}
la section nulle. Alors~:

\begin{proposition}\label{disque_1-convexe}
Avec ces notations, si $X$ est $1$-convexe et si la m\'etrique
continue $\|.\|$ sur $L$ est \`a courbure strictement positive,
le fibr\'e en disques $D_r(X,L)$ est lui aussi $1$-convexe. 
\end{proposition}
\begin{proof}
%Rappelons en effet que la stricte positivit\'e de la courbure de $\|.\|$
%signifie que pour toute section $l$ de $L$
%ne s'annulant pas sur un ouvert $U$ de $X$, la fonction $-\log\|l\|$
%est strictement \psh sur $U$. Pla\c{c}ons-nous sous les
%hypoth\`eses de la proposition, et supposons donc 
Supposons donn\'ee une fonction
d'exhaustion continue $\phi$ sur $X$ strictement \psh
hors d'un compact $K$, et
consid\'erons aussi une fonction strictement convexe
continue $\chi:[-\infty,\log r[\longto\R_+$
tendant vers $0$ en $-\infty$ et vers $+\infty$ en $\log r$.
Alors il d\'ecoule de l'hypoth\`ese de stricte positivit\'e sur la
courbure de $\|.\|$ que
la fonction $\psi$ d\'efinie pour $z\in D_r(X,L)$ par la formule
\begin{equation}
\psi(z)=\phi(\pi(z))+\chi(\log\|z\|^\vee),
\end{equation}
continue et exhaustive sur $D_r(X,L)$, est bien
strictement \psh hors du compact $\iota(K)$.
\end{proof}

Rappelons enfin la version suivante du th\'eor\`eme de
l'application ouverte de Banach (cf. p.~ex. \cite{BouEVT}, I~\S3,
cor.~3 p.~I.19 et ex.~4 p.~I.28)~:
\begin{theoreme}\label{application_ouverte}
Soient $E$ et $F$ deux espaces de Fr\'echet et $u:E\longto F$
une application lin\'eaire continue \`a conoyau
de dimension finie. Alors le sous-espace $u(E)$ de $F$ est
ferm\'e et l'application $u:E\longto u(E)$ est ouverte.
En particulier, pour toute semi-norme continue $p$ sur $E$
il existe une semi-norme continue $q$ sur $F$ telle que,
pour tout $y$ dans $u(E)$, il existe $x$ dans $E$ v\'erifiant
\begin{equation}
u(x)=y
\end{equation}
et
\begin{equation}
p(x)\leq q(y).
\end{equation}
\end{theoreme}

\'Enon\c{c}ons maintenant le r\'esultat principal de cette section~:
\begin{theoreme}
\label{comparaison_normes}
Soient $X$ un espace analytique $1$-convexe
r\'eduit, $\overline{L}$ un $\cO_X$-module inversible
hermitien \`a courbure strictement positive, $\overline{E}$
un $\cO_X$-module localement libre hermitien de type fini, $i:Y\inj X$
un sous-espace analytique ferm\'e r\'eduit de $X$,
et $\overline{V}$ un $\cO_Y$-module localement libre hermitien de type fini.
On suppose donn\'es un sous-$\cO_X$-module
coh\'erent $F$ de $E$ et un morphisme surjectif
\begin{equation}
p:F\surj i_*V
\end{equation}
de $\cO_X$-modules coh\'erents.
Alors~:
\begin{enumerate}
\item
Pour
tout r\'eel $\epsilon>0$ et pour tout compact non
vide $B$ de $Y$, il existe un r\'eel $C>0$ et un compact
non vide $A$ de $X$ tels que
pour tout $n\geq 0$, pour tout $s\in\Gamma(Y,V\tens L^{\tens n})$
et pour tout $\widetilde{s}\in\Gamma(X,F\tens L^{\tens n})$
v\'erifiant $p(\widetilde{s})=i_*s$,
on ait
\begin{equation}
\|\widetilde{s}\|_{L^\infty(A)}\geq C e^{-n\epsilon}\|s\|_{L^\infty(B)}.
\end{equation}
\item
Il existe un entier $n_0$ et, pour
tout r\'eel $\epsilon>0$ et pour tout compact non
vide $A$ de $X$, un r\'eel $C'>0$ et un compact
non vide $B$ de $Y$, tels que
pour tout $n\geq n_0$ et pour tout $s\in\Gamma(Y,V\tens L^{\tens n})$
il existe $\widetilde{s}\in\Gamma(X,F\tens L^{\tens n})$
v\'erifiant $p(\widetilde{s})=i_*s$
et
\begin{equation}
\|\widetilde{s}\|_{L^\infty(A)}\leq C'e^{n\epsilon}\|s\|_{L^\infty(B)}.
\end{equation}
\end{enumerate}
\end{theoreme}

La preuve de ce th\'eor\`eme,
directement adapt\'ee de \cite{Bost}, va occuper les paragraphes
\ref{debut_preuve} \`a \ref{fin_preuve}~:

\entry\label{debut_preuve}
Conservons les notations introduites en \ref{notations_disque_1-convexe},
en posant pour all\'eger l'\'ecriture
\begin{equation}
D(X)=D_1(X,L)
\end{equation}
et, si $A$ est une partie de $X$, $D_r(A,L)$ l'ensemble des
points de $D_r(X,L)$ dont l'image 
par la projection
canonique $\pi$ est dans $A$~;
on adoptera aussi les notations analogues sur $Y$.
Remarquons notamment que $D_r(Y,L)$
s'identifie au produit fibr\'e de $D_r(X,L)$ et
de $Y$ au-dessus de $X$. On notera encore
\begin{equation}
i:D(Y)\inj D(X)
\end{equation}
l'immersion ferm\'ee naturelle.

Par \ref{rappel_topologie_canonique}(i), la famille
de semi-normes
\begin{equation}
\|.\|_{D(X),A,r}=\|.\|_{L^\infty(D_r(A,L),\pi^*E)}
\end{equation}
(resp. $\|.\|_{D(Y),B,r}=\|.\|_{L^\infty(D_r(B,L),\pi^*V)}$),
pour $A$ compact non vide de $X$ 
(resp. $B$ compact non vide de $Y$) et $r\in]0,1[$,
d\'efinit la topologie de Fr\'echet canonique sur $\Gamma(D(X),\pi^*E)$
(resp. sur $\Gamma(D(Y),\pi^*V)$).

La projection $\pi:D(X)\surj X$ \'etant plate, l'injection $F\inj E$
sur $X$ se rel\`eve en $\pi^*F\inj\pi^*E$ sur $D(X)$ et,
passant aux sections globales, il r\'esulte
de \ref{rappel_topologie_canonique}(ii) que la topologie canonique sur
\begin{equation}
\label{inclusionD(X)FE}
\Gamma(D(X),\pi^*F)\inj\Gamma(D(X),\pi^*E)
\end{equation}
est encore d\'efinie par la famille des semi-normes $\|.\|_{D(X),A,r}$.

\entry %(action de U(1) et projections...)
Remarquons que l'action naturelle du groupe $\mathbb{G}_m$
sur $V(X,L)$, d\'efinie par l'action des homoth\'eties sur
les fibres de $L^\vee$, se restreint en une action continue du groupe
unitaire $U(1)=\{u\in\C\;|\;|u|=1\}$ sur $D(X)$.
Pour tout entier $k$ notons alors
\begin{equation}
\Gamma(D(X),\pi^*E)_k\subset\Gamma(D(X),\pi^*E)
\end{equation}
le sous-espace form\'e des sections analytiques $f$ de $\pi^*E$
sur $D(X)$ v\'erifiant
\begin{equation}
f(uz)=u^kf(z)
\end{equation}
pour tous $u\in U(1)$ et $z\in D(X)$, o\`u l'on a identifi\'e
les fibres $(\pi^*E)_{z}$ et $(\pi^*E)_{uz}$
\`a leur image commune $E_{\pi(z)}=E_{\pi(uz)}$ par $\pi$.

Notons aussi
\begin{equation}
q_k:\Gamma(D(X),\pi^*E)\surj\Gamma(D(X),\pi^*E)_k
\end{equation}
la projection d\'efinie par la formule
\begin{equation}
q_k(f)(z)=\int_0^1e^{-2i\pi kt}f(e^{2i\pi t}z)dt\;\in E_{\pi(z)}.
\end{equation}
Pour tous $A$ compact de $X$ et $r\in]0,1[$ on a clairement
\begin{equation}
\label{qk_contractante}
\|q_k(f)\|_{D(X),A,r}\leq\|f\|_{D(X),A,r},
\end{equation}
de sorte que $q_k$ est continue.

Remarquons enfin que le sous-espace $\Gamma(D(X),\pi^*E)_k$
s'identifie naturellement \`a l'espace $\Gamma(X,E\tens L^{\tens k})$
des sections analytiques de $E\tens L^{\tens k}$ sur $X$, cette
identification associant \`a l'\'el\'ement $s\in\Gamma(X,E\tens L^{\tens k})$
la section analytique $f$ de $\pi^*E$ sur $D(X)$ d\'efinie pour $z\in D(X)$
par
\begin{equation}
\label{fs}
f(z)=<s(\pi(z)),z^{\tens k}>\;\in E_{\pi(z)},
\end {equation}
o\`u $s(\pi(z))$ appartient \`a la fibre $(E\tens L^{\tens k})_{\pi(z)}$
et o\`u $z^{\tens k}$ est consid\'er\'e comme un \'el\'ement
de la droite complexe $L^{\vee\tens k}_{\pi(z)}$. 
Avec ces notations, on a alors clairement
\begin{equation}
\label{normes_fs}
\|f\|_{D(X),A,r}=r^k\|s\|_{L^\infty(A)}
\end{equation}
pour tous $A$ compact de $X$ et $r\in]0,1[$.

\entry
De la m\^eme fa\c{c}on, on dispose d'une 
action continue naturelle de $U(1)$
sur $\Gamma(D(X),\pi^*F)$ (resp. sur $\Gamma(D(Y),\pi^*V)$) \`a
laquelle sont encore associ\'es des sous-espaces
propres $\Gamma(D(X),\pi^*F)_k$ (resp. $\Gamma(D(Y),\pi^*V)_k$)
s'identifiant naturellement \`a $\Gamma(X,F\tens L^{\tens k})$
(resp. \`a $\Gamma(Y,V\tens L^{\tens k})$)
et des projections $q_k$ d\'efinies par les m\^emes formules,
de sorte que les in\'egalit\'es
et \'egalit\'es analogues \`a \eqref{qk_contractante}
et \eqref{normes_fs} restent valides.

Remarquons notamment que l'inclusion \eqref{inclusionD(X)FE}
est compatible \`a toutes ces constructions.

\entry
En utilisant \`a nouveau la platitude de $\pi$ pour
relever sur $D(X)$ la suite exacte
courte $0\longto\ker p\longto F\longto i_*V\longto 0$, et en
passant \`a la cohomologie, on obtient une suite exacte
\begin{equation}
\Gamma(D(X),\pi^*F)\longto\Gamma(D(X),\pi^*i_*V)\longto H^1(D(X),\pi^*\ker p)
\end{equation}
o\`u, par la proposition \ref{disque_1-convexe} et par la
caract\'erisation \ref{def_1-convexe}(iv) des espaces $1$-convexes,
le groupe de cohomologie $H^1(D(X),\pi^*\ker p)$ est de dimension finie.
Compte tenu de \ref{rappel_topologie_canonique}(iii),
{}$\Gamma(D(X),\pi^*i_*V)=\Gamma(D(X),i_*\pi^*V)$ s'identifie
hom\'eomorphiquement \`a $\Gamma(D(Y),\pi^*V)$, et l'on en d\'eduit
une application lin\'eaire continue 
\begin{equation}
\rho:\Gamma(D(X),\pi^*F)\longto\Gamma(D(Y),\pi^*V)
\end{equation}
\`a conoyau de dimension finie.

\entry
La continuit\'e de $\rho$ implique que, 
pour tout r\'eel $\epsilon>0$ et pour tout compact $B$ de $Y$,
il existe $C>0$, $r\in]0,1[$ et $A$ compact de $X$ tels que,
pour tout $\widetilde{f}\in\Gamma(D(X),\pi^*F)$,
posant $f=\rho(\widetilde{f})$, on ait
\begin{equation}
\label{f<ftilde} 
\|f\|_{D(Y),B,e^{-\epsilon}}\leq C^{-1}\|\widetilde{f}\|_{D(X),A,r}.
\end{equation}

Soient maintenant $n\geq0$, $s\in\Gamma(Y,V\tens L^{\tens n})$
et $\widetilde{s}\in\Gamma(X,F\tens L^{\tens n})$
v\'erifiant $p(\widetilde{s})=i_*s$.
Si l'on note $f\in\Gamma(D(Y),\pi^*V)_n$
et $\widetilde{f}\in\Gamma(D(X),\pi^*F)_n$ les \'el\'ements
qui leur sont associ\'es par l'analogue de \eqref{fs}, la commutativit\'e
du diagramme
\begin{equation}
\label{diagramme_commutatif}
\begin{CD}
\Gamma(D(X),\pi^*F)_n @>\rho>> \Gamma(D(Y),\pi^*V)_n\\
@| @|\\
\Gamma(X,F\tens L^{\tens n}) @>{(i_*)^{-1}\circ p}>> \Gamma(Y,V\tens L^{\tens n})
\end{CD}
\end{equation}
implique qu'on a bien $f=\rho(\widetilde{f})$,
et \eqref{normes_fs}
et \eqref{f<ftilde} donnent
\begin{equation}
\begin{split}
\|\widetilde{s}\|_{L^\infty(A)}&=\|\widetilde{f}\|_{D(X),A,1}\\
&\geq\|\widetilde{f}\|_{D(X),A,r}\\
&\geq C \|f\|_{D(Y),B,e^{-\epsilon}}\\
&=C e^{-n\epsilon}\|s\|_{L^\infty(B)},
\end{split}
\end{equation}
ce qui d\'emontre le point (i) du th\'eor\`eme.

\entry
D'autre part, on peut aussi appliquer le
th\'eor\`eme \ref{application_ouverte}
\`a l'application $\rho$~: on trouve ainsi que
l'image $\im(\rho)=\rho(\Gamma(D(X),\pi^*F))$
de $\rho$ dans $\Gamma(D(Y),\pi^*V)$
est ferm\'ee (de sorte que la topologie quotient sur $\coker\rho$
co{\"\i}ncide avec sa topologie usuelle de $\C$-espace vectoriel
de dimension finie) et que, pour tout $\epsilon>0$ et pour tout
compact $A$ de $X$,
il existe $C'>0$, $r'\in]0,1[$ et $B$ compact de $Y$ tels que,
pour tout $f\in\im(\rho)$, il existe $\widetilde{f}\in\Gamma(D(X),\pi^*F)$
v\'erifiant $f=\rho(\widetilde{f})$ et
\begin{equation}
\label{ftilde<f}
\|\widetilde{f}\|_{D(X),A,e^{-\epsilon}}\leq C'\|f\|_{D(Y),B,r'}.
\end{equation}

Par ailleurs, l'application $\rho$ \'etant clairement $U(1)$-\'equivariante,
on en d\'eduit par passage au quotient une action continue de $U(1)$
sur $\coker\rho$.
Ceci permet de consid\'erer la d\'ecomposition
\begin{equation}
\coker\rho=\bigoplus_n(\coker\rho)_n
\end{equation}
de ce conoyau en somme (finie~!) de sous-espaces propres
pour cette action, o\`u $(\coker\rho)_n$ est l'espace
propre associ\'e au caract\`ere $(u\mapsto u^n)$ de $U(1)$.
Notons alors $n_0$ le plus petit entier tel que $(\coker\rho)_n$
soit nul pour $n\geq n_0$.

\entry
Donnons-nous maintenant un r\'eel $\epsilon>0$,
un compact $A$ de $X$,
un entier $n\geq n_0$,
et un \'el\'ement $s\in\Gamma(Y,V\tens L^{\tens n})$
auquel par l'analogue de \eqref{fs} 
on peut associer $f\in\Gamma(D(Y),\pi^*V)_n$.
Puisque $(\coker\rho)_n$ est nul, $f$ appartient
\`a $\im(\rho)$, et par la discussion qui
pr\'ec\`ede on peut \'ecrire $f=\rho(\widetilde{f})$
o\`u $\widetilde{f}\in\Gamma(D(X),\pi^*F)$ v\'erifie \eqref{ftilde<f}.
L'application $\rho$ \'etant compatible aux projections $q_n$,
on a encore $f=\rho(q_n(\widetilde{f}))$ et,
le diagramme \eqref{diagramme_commutatif}
\'etant commutatif, notant $\widetilde{s}\in\Gamma(X,F\tens L^{\tens n})$
l'\'el\'ement associ\'e \`a $q_n(\widetilde{f})$, on trouve
\begin{equation}
p(\widetilde{s})=i_*s
\end{equation}
avec, par \eqref{qk_contractante}, \eqref{normes_fs}
et \eqref{ftilde<f},
\begin{equation}
\begin{split}
\|\widetilde{s}\|_{L^\infty(A)}&=e^{n\epsilon}\|q_n(\widetilde{f})\|_{D(X),A,e^{-\epsilon}}\\
&\leq e^{n\epsilon}\|\widetilde{f}\|_{D(X),A.e^{-\epsilon}}\\
&\leq C' e^{n\epsilon}\|f\|_{D(Y),B,r'}\\
&\leq C' e^{n\epsilon}\|f\|_{D(Y),B,1}\\
&=C' e^{n\epsilon}\|s\|_{L^\infty(B)},
\end{split}
\end{equation}
ce qui d\'emontre le point (ii).
\label{fin_preuve}

La preuve du th\'eor\`eme \'etant achev\'ee, indiquons maintenant
comment celui-ci
peut s'interpr\'eter en termes de normes sous-quotient.

\entry
\label{hypotheses_corollaire}
Soient $X$ un espace analytique \emph{compact}
r\'eduit, $\overline{L}$ un $\cO_X$-module inversible
hermitien \`a courbure strictement positive, $\overline{E}$
un $\cO_X$-module localement libre hermitien de type fini, $i:Y\inj X$
un sous-espace analytique ferm\'e r\'eduit de $X$,
et $\overline{V}$ un $\cO_Y$-module localement libre hermitien de type fini.

On suppose le $\cO_X$-module coh\'erent $i_*V$ muni d'une structure
de sous-quotient de $E$~; rappelons que ceci correspond \`a la
donn\'ee d'un isomorphisme $\phi:i_*V\simeq E_1/E_2$
o\`u $E\supset E_1\supset E_2$ est une filtration \`a deux termes
de $E$ par des sous-$\cO_X$-modules coh\'erents.

Le faisceau $L$ \'etant ample,
le morphisme naturel
\begin{equation}
\psi:\Gamma(X,E_1\tens L^{\tens n})/\Gamma(X,E_2\tens L^{\tens n})\longto\Gamma(X,(E_1/E_2)\tens L^{\tens n})
\end{equation}
est inversible pour tout $n$ assez grand, de sorte que le
morphisme compos\'e
\begin{equation}
\psi^{-1}\circ\phi:\Gamma(X,i_*V\tens L^{\tens n})\overset{\sim}{\longto}\Gamma(X,E_1\tens L^{\tens n})/\Gamma(X,E_2\tens\cL^{\tens n})
\end{equation}
munit $\Gamma(X,i_*V\tens L^{\tens n})$ d'une structure de
sous-quotient de $\Gamma(X,E\tens L^{\tens n})$, relativement
\`a la filtration 
\begin{equation}
\Gamma(X,E\tens L^{\tens n})\supset\Gamma(X,E_1\tens L^{\tens n})\supset\Gamma(X,E_2\tens L^{\tens n}).
\end{equation}
Ainsi le $\C$-espace vectoriel $\Gamma(X,i_*V\tens L^{\tens n})=\Gamma(Y,V\tens i^*L^{\tens n})$ dispose-t-il~:
\begin{itemize}
\item de la norme uniforme $\|.\|_{L^\infty(Y)}$ sur $Y$
\item de la norme $\|.\|_{\mathrm{sq},L^\infty(X)}$ provenant par la
structure de sous-quotient de
la norme uniforme sur $X$.
\end{itemize}
Alors~:

\begin{corollaire}
\label{cor_presque-isom}
Sous les hypoth\`eses \ref{hypotheses_corollaire},
il existe un entier $n_0$ et, pour
tout r\'eel $\epsilon>0$, deux r\'eels $C>0$ et $C'>0$, tels que
pour tout $n\geq n_0$ les deux normes sur $\Gamma(X,i_*V\tens L^{\tens n})$
ainsi construites
v\'erifient les in\'egalit\'es 
\begin{equation}
\label{Y<sqX}
\|.\|_{\mathrm{sq},L^\infty(X)}\geq Ce^{-n\epsilon}\|.\|_{L^\infty(Y)}
\end{equation}
et
\begin{equation}
\label{sqX<Y}
\|.\|_{\mathrm{sq},L^\infty(X)}\leq C'e^{n\epsilon}\|.\|_{L^\infty(Y)}.
\end{equation}
\end{corollaire}
\begin{proof}
Notant $F=E_1$ et $p:F\surj i_*V$ la projection d\'efinissant
la structure de sous-quotient sur $i_*V$
(suivant la caract\'erisation \ref{caracterisations_sous-quotients}(i)),
on se retrouve en position d'appliquer le th\'eor\`eme.
Compte tenu de la d\'efinition \eqref{formule_norme} des normes
sous-quotient, l'in\'egalit\'e \eqref{Y<sqX} n'est autre qu'une
traduction du point (i) du th\'eor\`eme o\`u l'on a choisi $B=Y$,
qui est bien compact puisque $X$ l'est par hypoth\`ese.
De la m\^eme fa\c{c}on, l'in\'egalit\'e \eqref{sqX<Y}
r\'esulte du point (ii) du th\'eor\`eme avec $A=X$.
\end{proof}

\section{Th\'eor\`eme de Hilbert-Samuel arithm\'etique}

\entry
On supposera ici le lecteur familier avec le langage de
la g\'eom\'etrie d'Arakelov.
N\'eanmoins, afin d'\'eviter toute ambigu\"{\i}t\'e, on
commence par rappeler certaines notions de base et
pr\'eciser les normalisations qui seront utilis\'ees.

Soit $K$ un corps de nombres d'anneau d'entiers $\cO_K$.
Un $\cO_K$-module norm\'e $\overline{M}$ est
la donn\'ee d'un $\cO_K$-module de type fini $M$ et, pour
tout plongement $\sigma$ de $K$ dans $\C$, d'une norme $\|.\|_\sigma$
sur le $\C$-espace vectoriel $M_\sigma=M\tens_{\cO_K,\sigma}\C$,
ceci de fa\c{c}on compatible \`a la conjugaison complexe.
Notons alors $M_{\textrm{tors}}$ le sous-module de torsion
de $M$, $M_{\textrm{libre}}=M/M_{\textrm{tors}}$ son plus grand
quotient sans torsion, qui s'identifie \`a un r\'eseau
de $M_\R=M\tens_\Z\R$, et enfin $B\subset M_\R$ la boule unit\'e
de $M_\R$, c'est-\`a-dire l'ensemble des
\'el\'ements $m\in M_\R$ dont les images par les applications
naturelles $M_\R\longto M_\sigma$ sont de norme $\|m\|_\sigma\leq1$
pour tout $\sigma$. Avec ces notations, le degr\'e d'Arakelov
de $\overline{M}$ est le r\'eel
\begin{equation}
\wdeg\overline{M}=\log\# M_{\textrm{tors}}-\log\vol(M_\R/M_{\textrm{libre}})
\end{equation}
o\`u le volume $\vol(M_\R/M_{\textrm{libre}})$ est pris relativement
\`a l'unique mesure de Haar sur $M_\R$ qui donne \`a $B$ le volume $1$.
Il pourra aussi s'av\'erer commode de consid\'erer le degr\'e d'Arakelov
normalis\'e $\frac{1}{[K:\Q]}\wdeg$,
comme on le trouve parfois dans la litt\'erature.

\entry
On v\'erifie ais\'ement que si $(\|.\|_\sigma)_\sigma$
et $(\|.\|'_\sigma)_\sigma$ sont deux familles de normes munissant
un m\^eme $\cO_K$-module de type fini $M$ de deux structures
de $\cO_K$-module norm\'e $\overline{M}$ et $\overline{M}'$,
et si $C>0$ est un r\'eel tel que pour tout $\sigma$
on ait
\begin{equation}
\|.\|_\sigma\leq C\|.\|'_\sigma,
\end{equation}
alors
\begin{equation}
\label{inegalite_wdeg}
\wdeg\overline{M}\geq\wdeg\overline{M}'-\rg_\Z M\log C
\end{equation}
o\`u $\rg_\Z M=\dim_\R M_\R=[K:\Q]\rg_{\cO_K}M$.

\entry
Lorsque les normes $\|.\|_\sigma$ proviennent de produits
scalaires hermitiens $(.,.)_\sigma$ sur $M_\sigma$, on dit que
le $\cO_K$-module norm\'e $\overline{M}$ est
un $\cO_K$-module hermitien. On v\'erifie alors que le degr\'e
d'Arakelov de $\overline{M}$ peut aussi s'exprimer sous la forme
\begin{equation}
\wdeg\overline{M}=\log\# M/(s_1,\dots,s_r)-\sum_\sigma\log\|s_1\wedge\dots\wedge s_r\|_{\bigwedge^r M_\sigma}
\end{equation}
o\`u $s_1,\dots,s_r$ sont des éléments de $M$ dont
les images dans le $K$-espace vectoriel $M_K$ en forment une base,
et o\`u $\|.\|_{\bigwedge^r M_\sigma}$ est la norme
hermitienne puissance ext\'erieure $r$-i\`eme de $\|.\|_\sigma$.

Si $\overline{M}$ est un $\cO_K$-module hermitien, et si
\begin{equation}
0\longto N\longto M\longto Q\longto 0
\end{equation}
est une suite exacte courte de $\cO_K$-modules de type fini,
les normes restreintes (resp. quotient) sur $N$ (resp. sur $Q$)
sont encore hermitiennes, et pour les structures de $\cO_K$-modules
hermitiens correspondantes on a
\begin{equation}
\label{additivite}
\wdeg\overline{M}=\wdeg\overline{N}+\wdeg\overline{Q}.
\end{equation}

\entry
\label{notations_HSA}
Consid\'erons maintenant $\gX$ un
sch\'ema projectif sur $\spec\cO_K$
de fibre g\'en\'erique $\gX_K$ r\'eduite, $\overline{\cL}$ 
un $\cO_{\gX}$-module inversible ample hermitien \`a
courbure semi-positive, et $\overline{\cC}$
un $\cO_{\gX}$-module coh\'erent muni d'une structure de
sous-quotient d'un $\cO_{\gX}$-module localement libre hermitien~;
plus pr\'ecis\'ement, $\overline{\cC}$ correspond \`a la donn\'ee~:
\begin{itemize}
\item d'un $\cO_{\gX}$-module coh\'erent $\cC$,
\item d'un $\cO_{\gX}$-module localement libre hermitien $\overline{\cE}$,
\item d'une filtration \`a deux termes $\cE\supset\cE_1\supset\cE_2$
de $\cE$ par des sous-$\cO_{\gX}$-modules coh\'erents, et
\item d'un isomorphisme $\phi:\cC\overset{\sim}{\longto}\cE_1/\cE_2$.
\end{itemize}
On dira que $\overline{\cC}$ est un $\cO_{\gX}$-module coh\'erent muni
de m\'etriques de sous-quotient.

Notons $d=\dim|\cC|$ la dimension du support de $\cC$
(ainsi $d\leq\dim\gX$) et
\begin{equation}
[\cC]=\sum_{x\in\gX_{(d)}}(\lg_{\cO_{\gX,x}}\cC_x)[x]
\end{equation}
le cycle de dimension $d$ associ\'e. On supposera enfin $\gX(\C)$
muni d'une forme continue $\mu$ d\'efinissant une forme
volume strictement positive en tout point du lieu lisse de $\gX(\C)$
et compatible \`a la conjugaison complexe
(remarquons qu'une telle forme existe sur tout
espace analytique complexe projectif r\'eduit). On aura alors
\`a consid\'erer l'une des deux situations suivantes~:
\begin{itemize}
\item[\emph{Hypoth\`ese} ($L^\infty$)~:]
pour tout entier $n$, $\Gamma(\gX,\cE\tens\cL^{\tens n})$ est muni d'une
structure de $\cO_K$-module norm\'e au moyen des normes uniformes
sur les $\gX_\sigma(\C)$.
\item[\emph{Hypoth\`ese} ($L^2(\mu)$)~:]
pour tout entier $n$, $\Gamma(\gX,\cE\tens\cL^{\tens n})$
est muni d'une
structure de $\cO_K$-module hermitien au moyen des normes
d'int\'egration $L^2(\mu)$
sur les $\gX_\sigma(\C)$.
\end{itemize}
Munissons alors $\Gamma(\gX,\cE_1\tens\cL^{\tens n})$
et $\Gamma(\gX,\cE_2\tens\cL^{\tens n})$ des normes restreintes,
et $\Gamma(\gX,\cE_1\tens\cL^{\tens n})/\Gamma(\gX,\cE_2\tens\cL^{\tens n})$
des normes quotient.

\begin{definition}
\label{def_HSA}
La fonction de Hilbert-Samuel arithm\'etique
\begin{equation}
h(\overline{\cC};.):\Z\longto\R
\end{equation}
du $\cO_{\gX}$-module coh\'erent muni de m\'etriques
de sous-quotient $\overline{\cC}$ est donn\'ee par la formule
\begin{equation}
h(\overline{\cC};n)=\wdeg\overline{\Gamma(\gX,\cE_1\tens\cL^{\tens n})/\Gamma(\gX,\cE_2\tens\cL^{\tens n})}.
\end{equation}
\end{definition}
Si l'on veut pr\'eciser sous laquelle des hypoth\`eses ($L^\infty$) ou
($L^2(\mu)$) on s'est plac\'e pour construire cette fonction de
Hilbert-Samuel arithm\'etique, on pourra noter
celle-ci $h_{L^\infty}(\overline{\cC};.)$
ou $h_{L^2(\mu)}(\overline{\cC};.)$, respectivement.

\entry
Remarquons que, les normes $L^2(\mu)$ \'etant hermitiennes,
il r\'esulte de \eqref{additivite} qu'on peut aussi \'ecrire
\begin{equation}
h_{L^2(\mu)}(\overline{\cC};n)=\wdeg\overline{\Gamma(\gX,\cE_1\tens\cL^{\tens n})}_{L^2(\mu)}-\wdeg\overline{\Gamma(\gX,\cE_2\tens\cL^{\tens n})}_{L^2(\mu)}.
\end{equation}
L'\'egalit\'e analogue pour les normes $L^\infty$
n'a en revanche a priori pas de raison d'\^etre v\'erifi\'ee.

\entry
Dans tous les cas, remarquons aussi que, le faisceau $\cL$
\'etant ample, le morphisme naturel
\begin{equation}
\psi:\Gamma(\gX,\cE_1\tens\cL^{\tens n})/\Gamma(\gX,\cE_2\tens\cL^{\tens n})\longto\Gamma(\gX,(\cE_1/\cE_2)\tens\cL^{\tens n})
\end{equation}
est inversible pour tout $n$ assez grand, de sorte que le
morphisme compos\'e
\begin{equation}
\psi^{-1}\circ\phi:\Gamma(\gX,\cC\tens\cL^{\tens n})\overset{\sim}{\longto}\Gamma(\gX,\cE_1\tens\cL^{\tens n})/\Gamma(\gX,\cE_2\tens\cL^{\tens n})
\end{equation}
munit $\Gamma(\gX,\cC\tens\cL^{\tens n})$ d'une structure de
sous-quotient de $\Gamma(\gX,\cE\tens\cL^{\tens n})$, relativement
\`a la filtration
\begin{equation}
\Gamma(\gX,\cE\tens\cL^{\tens n})\supset\Gamma(\gX,\cE_1\tens\cL^{\tens n})\supset\Gamma(\gX,\cE_2\tens\cL^{\tens n}).
\end{equation}
Munissant alors $\Gamma(\gX,\cC\tens\cL^{\tens n})$ d'une structure
de $\cO_K$-module norm\'e au moyen des normes sous-quotient des normes
sur $\Gamma(\gX,\cE\tens\cL^{\tens n})$, on trouve~:
\begin{equation}
\label{autre_expression_HSA}
h(\overline{\cC};n)=\wdeg\overline{\Gamma(\gX,\cC\tens\cL^{\tens n})}.
\end{equation}

La suite de l'expos\'e fera appel aux deux r\'esultats classiques suivants.

\begin{lemme}
\label{approximation_metriques}
Soient $X$ un espace analytique complexe projectif r\'eduit
et $L$ un $\cO_X$-module inversible ample. Alors toute m\'etrique
continue \`a courbure semi-positive sur $L$ est limite uniforme
d\'ecroissante de m\'etriques $C^\infty$ \`a courbure strictement positive.
\end{lemme}
Il s'agit du th\'eor\`eme 4.6.1 de \cite{Maillot}, o\`u $X$ est
suppos\'e lisse.
Remarquons toutefois que la preuve donn\'ee dans \cite{Maillot}
n'utilise cette hypoth\`ese de lissit\'e
que pour montrer que toute fonction \psh continue sur $X$
est localement limite uniforme d\'ecroissante de fonctions
plu\-ri\-sous\-har\-mo\-niques $C^\infty$. Or, puisqu'une
fonction sur un espace analytique complexe
est \psh si et seulement si elle peut
s'\'ecrire localement comme restriction
d'une fonction \psh sur un espace lisse, ce r\'esultat
reste vrai en g\'en\'eral.

\begin{proposition}
Soient $X$ un espace analytique complexe projectif
r\'eduit muni d'une forme continue $\mu$ d\'efinissant
une forme volume strictement positivesur son lieu lisse, $\overline{L}$
un $\cO_X$-module inversible ample
muni d'une m\'etrique (continue) \`a courbure semi-positive,
et $\overline{E}$ un $\cO_X$-module localement libre hermitien.
Alors, pour tout $\epsilon>0$, il existe des constantes $c_\epsilon>0$
et $C_\epsilon>0$ telles que pour tout $n\geq0$ et pour
tout $s\in\Gamma(X,E\tens L^{\tens n})$, on ait
\begin{equation}
c_\epsilon e^{-n\epsilon}\|s\|_{L^\infty(X)}\leq\|s\|_{L^2(\mu,X)}\leq C_\epsilon e^{n\epsilon}\|s\|_{L^\infty(X)}.
\end{equation}
\end{proposition}
Lorsque $X$ est lisse, que les m\'etriques et la forme volume sont $C^\infty$,
et que la m\'etrique sur $L$ est \`a courbure strictement positive,
il s'agit une cons\'equence
directe du lemme 30 de \cite{GSRR}.
Remarquons que la preuve donn\'ee dans \cite{GSRR} reste encore valide
si l'on suppose que la forme volume s'annule \`a l'ordre fini le
long d'un ferm\'e analytique de codimension strictement positive,
et si la m\'etrique sur $E$ est continue.
Le cas g\'en\'eral s'en d\'eduit alors, en approchant la forme volume
continue par une forme volume $C^\infty$, en passant \`a
une r\'esolution des singularit\'es, et en approchant la
m\'etrique sur $L$ par une m\'etrique $C^\infty$ \`a courbure
strictement positive au moyen du lemme \ref{approximation_metriques}.

\begin{corollaire}
\label{hL2=hLinfty}
Avec les notations de \ref{notations_HSA} et \ref{def_HSA}, on a
\begin{equation}
|h_{L^\infty}(\overline{\cC};n)-h_{L^2(\mu)}(\overline{\cC};n)|=o(n^d)
\end{equation}
quand $n$ tend vers l'infini.
\end{corollaire}
\begin{proof}
Cela r\'esulte de la proposition \ref{compare_normes_sq},
de \eqref{inegalite_wdeg}, de \eqref{autre_expression_HSA},
de la proposition pr\'ec\'edente,
et du fait que d'apr\`es le th\'eor\`eme de Hilbert-Samuel g\'eom\'etrique,
le rang du $\cO_K$-module $\Gamma(\gX,\cC\tens\cL^{\tens n})$
est un polyn\^ome en $n$ de degr\'e au plus $d-1$.
\end{proof}

%Rappelons aussi le lemme de d\'evissage bien connu qui suit.

\begin{lemme}
\label{devissage}
Soient $X$ un sch\'ema noeth\'erien
quasi-projectif\footnote{\cad admettant un faisceau inversible ample}
et $\cC$ un $\cO_X$-module coh\'erent. Alors $\cC$ admet une filtration
d\'ecroissante finie
\begin{equation}
\cC=\cC_0\supset\cC_1\supset\dots\supset\cC_N=0
\end{equation}
telle que, pour tout $j\in\{1,\dots,N\}$,
il existe un sous-sch\'ema ferm\'e int\`egre $Z_j$ 
de $X$, un $\cO_{Z_j}$-module inversible $\cM_j$ et un isomorphisme
\begin{equation}
\cC_{j-1}/\cC_j\simeq i_{Z_j*}\cM_j
\end{equation}
o\`u $i_{Z_j}$ est le morphisme d'immersion de $Z_j$ dans $X$.

De fa\c{c}on plus pr\'ecise, si
un $\cO_X$-module inversible ample $\cL$ est donn\'e, on
peut prendre les $\cM_j$ de la forme $i_{Z_j}^*\cL^{\tens n_j}$
pour $n_j\in\Z$.
\end{lemme}
\begin{proof}
Proc\'edant par induction noeth\'erienne, on va montrer que
si $\cC'$ est un sous-$\cO_X$-module coh\'erent de $\cC$
maximal parmi ceux qui admettent une filtration v\'erifiant
les propri\'et\'es \'enonc\'ees dans le lemme, alors $\cC'=\cC$.
Dans le cas contraire, le $\cO_X$-module coh\'erent $\cC/\cC'$
est non nul donc, par \cite{EGA}~IV~3.1.5, poss\`ede un cycle
premier associ\'e $Z$~; autrement dit, il existe un voisinage $U$
du point g\'en\'erique $z$ de $Z$ dans $X$
et une section $s$ de $\cC/\cC'$ sur $U$ de support $Z\cap U$
(\cite{EGA}~IV~3.1.3(b)). D'autre part, par la caract\'erisation
\cite{EGA}~II~4.5.2(a) des faisceaux amples, on peut
apr\`es restriction supposer $U=X_f$ 
o\`u $f$ est une section globale de $\cL^{\tens n}$ sur $X$,
pour un certain entier $n$.
Alors, par \cite{EGA}~I~9.3.1(ii), il existe
un entier $m$ tel que $s\tens f^{\tens m}$ se
prolonge en une section globale de $(\cC/\cC')\tens\cL^{\tens mn}$
sur $X$, de sorte que $s\tens f^{\tens m+1}$ d\'efinit
une section globale de $(\cC/\cC')\tens\cL^{\tens (m+1)n}$
de support $Z$, ou encore une injection
\begin{equation}
\label{iZ*MinjC/C'}
i_{Z*}\cM\inj\cC/\cC',
\end{equation}
o\`u l'on a pos\'e $\cM=i_Z^*\cL^{\tens-(m+1)n}$.
Ainsi le sous-$\cO_X$-module de $\cC$ image inverse pour la projection
canonique de l'image de \eqref{iZ*MinjC/C'} dans $\cC/\cC'$
contient $\cC'$ et admet une filtration du type consid\'er\'e,
ce qui contredit l'hypoth\`ese de maximalit\'e faite sur $\cC'$.
\end{proof}

On est maintenant en mesure de prouver
le <<th\'eor\`eme de Hilbert-Samuel arithm\'etique>> annonc\'e.

\begin{theoreme}
Avec les notations de \ref{notations_HSA} et \ref{def_HSA}, on a
pour $n$ tendant vers l'infini le d\'eveloppement asymptotique
\begin{equation}
h(\overline{\cC};n)=\frac{n^d}{d!}(\wc1(\overline{\cL})^d.[\cC])+o(n^d),
\end{equation}
o\`u le nombre d'intersection arithm\'etique
g\'en\'eralis\'e $(\wc1(\overline{\cL})^d.[\cC])$ est d\'efini
au sens de \cite{Maillot}, \S5.
\end{theoreme}
\begin{proof}
Remarquons d'abord que par le corollaire \ref{hL2=hLinfty}, il est
indiff\'erent pour prouver le th\'eor\`eme de consid\'erer
les normes $L^\infty$ ou $L^2(\mu)$.
Compte tenu du lemme \ref{approximation_metriques}, on ne
perdra pas de g\'en\'eralit\'e en supposant que la m\'etrique
sur $\cL$ est $C^\infty$ \`a courbure strictement positive.

Le lemme \ref{devissage} nous donne une filtration finie
\begin{equation}
\label{tmp_filtr_C}
\cC=\cC_0\supset\cC_1\supset\dots\supset\cC_N=0
\end{equation}
de $\cC$ et, pour tout $j\in\{1,\dots,N\}$,
un sous-sch\'ema ferm\'e int\`egre $i_{Z_j}:Z_j\inj\gX$ 
et un $\cO_{Z_j}$-module inversible $\cM_j$ tel que $i_{Z_j*}\cM_j$
soit muni d'une structure de sous-quotient de $\cC$
\begin{equation}
i_{Z_j*}\cM_j\simeq\cC_{j-1}/\cC_j 
\end{equation}
et donc, par composition, d'une structure de sous-quotient de $\cE$
\begin{equation}
i_{Z_j*}\cM_j\simeq\widehat{\cC_{j-1}}/\widehat{\cC_j} 
\end{equation}
o\`u
\begin{equation}
\cE_1=\widehat{\cC_0}\supset\widehat{\cC_1}\supset\dots\supset\widehat{\cC_N}=\cE_2
\end{equation}
est la filtration relevant \eqref{tmp_filtr_C} via
l'isomorphisme $\phi:\cC\simeq\cE_1/\cE_2$.
On notera $\overline{i_{Z_j*}\cM_j}$
le $\cO_{\gX}$-module coh\'erent muni de m\'etriques de sous-quotient
ainsi obtenu.
%Soit $n$ un entier assez grand pour que les morphismes
%naturels $\Gamma(\gX,\cE_1\tens\cL^{\tens n})\longto\Gamma(\gX,(\cE_1/\cE_2)\tens\cL^{\tens n})$
%et $\Gamma(\gX,\cC_{j-1}\tens\cL^{\tens n})\longto\Gamma(\gX,(\cC_{j-1}/\cC_j)\tens\cL^{\tens n})$
%soient surjectifs.
Par la proposition \ref{transitivite_norme}, les normes
sur $\Gamma(\gX,\widehat{\cC_{j-1}}\tens\cL^{\tens n})/\Gamma(\gX,\widehat{\cC_j}\tens\cL^{\tens n})$
obtenues lorsqu'on le consid\`ere comme sous-quotient
de $\Gamma(\gX,\cE_1\tens\cL^{\tens n})/\Gamma(\gX,\cE_2\tens\cL^{\tens n})$
et celles obtenues lorsqu'on le consid\`ere
comme sous-quotient de $\Gamma(\gX,\cE\tens\cL^{\tens n})$
co\"{\i}ncident~; une application r\'ep\'et\'ee de \eqref{additivite}
donne alors
\begin{equation}
\begin{split}
h_{L^2(\mu)}(\overline{\cC};n)&=\wdeg\overline{\Gamma(\gX,\cE_1\tens\cL^{\tens n})/\Gamma(\gX,\cE_2\tens\cL^{\tens n})}_{L^2(\mu)}\\
&=\wdeg\overline{\Gamma(\gX,\widehat{\cC_0}\tens\cL^{\tens n})/\Gamma(\gX,\widehat{\cC_1}\tens\cL^{\tens n})}_{L^2(\mu)}+\dots\\
&\phantom{=========}+\wdeg\overline{\Gamma(\gX,\widehat{\cC_{N-1}}\tens\cL^{\tens n})/\Gamma(\gX,\widehat{\cC_N}\tens\cL^{\tens n})}_{L^2(\mu)}\\
&=h_{L^2(\mu)}(\overline{i_{Z_1*}\cM_1};n)+\dots+h_{L^2(\mu)}(\overline{i_{Z_N*}\cM_N};n).
\end{split}
\end{equation}
Les supports des $i_{Z_j*}\cM_j$ \'etant tous de dimension au plus $d$,
ceci implique que le th\'eor\`eme sera vrai pour $\overline{\cC}$
d\`es lors qu'il sera vrai pour les $\overline{i_{Z_j*}\cM_j}$.
Ainsi on pourra maintenant supposer $\cC=i_{Z*}\cM$
o\`u $i:Z\inj\gX$ est un sous-sch\'ema ferm\'e int\`egre
de dimension $d$
et $\cM$ un $\cO_Z$-module inversible. 

Sous cette hypoth\`ese, consid\'erons un entier $n$ suffisamment grand
pour que le morphisme naturel
\begin{equation}
\Gamma(\gX,\cE_1\tens\cL^{\tens n})\longto\Gamma(\gX,(\cE_1/\cE_2)\tens\cL^{\tens n})\simeq\Gamma(\gX,(i_*\cM)\tens\cL^{\tens n})
\end{equation}
soit surjectif.
Le $\cO_K$-module $\Gamma(\gX,(i_*\cM)\tens\cL^{\tens n})=\Gamma(Z,\cM\tens i^*\cL^{\tens n})$ dispose alors
\begin{itemize}
\item des normes uniformes $\|.\|_{L^\infty(Z),\sigma}$ sur $Z$
\item des normes $\|.\|_{\mathrm{sq},L^\infty(\gX),\sigma}$ provenant
par la structure de sous-quotient des normes uniformes sur $\gX$.
\end{itemize}
Le corollaire \ref{cor_presque-isom}
permet de comparer ces deux normes, et
les m\^emes arguments que ceux de la preuve du corollaire \ref{hL2=hLinfty}
donnent pour $n$ tendant vers l'infini l'estimation
\begin{equation}
|\wdeg\overline{\Gamma(Z,\cM\tens i^*\cL^{\tens n})}_{\mathrm{sq},L^\infty(\gX)}-\wdeg\overline{\Gamma(Z,\cM\tens i^*\cL^{\tens n})}_{L^\infty(Z)}|=o(n^d).
\end{equation}
Puisque par d\'efinition
on a $\wdeg\overline{\Gamma(Z,\cM\tens i^*\cL^{\tens n})}_{\mathrm{sq},L^\infty(\gX)}=h_{L^\infty}(\overline{\cC};n)$,
il suffit pour conclure de v\'erifier
\begin{equation}
\wdeg\overline{\Gamma(Z,\cM\tens i^*\cL^{\tens n})}_{L^\infty(Z)}=\frac{n^d}{d!}(\wc1(\overline{\cL})^d.Z)+o(n^d).
\end{equation}
Or ceci n'est autre que le th\'eor\`eme 1.4 de \cite{Zhang}.
\end{proof}

\end{document}